\documentclass[11pt]{article}
\usepackage{amsmath,amssymb}
\usepackage{epsf}
\input epsf
\newcommand{\proof}{{\noindent \bf Proof. }}
\newtheorem{thm}{Theorem}[section]
\newtheorem{lem}[thm]{Lemma}

\newtheorem{quest}[thm]{Question}

\newcommand{\ve}{\varepsilon} 

\newif\ifnotesw\noteswtrue

\def\2{\mathbb Z_2}

\def\qed{\hfill$\Box$}
\def\deg{{\mathrm deg}}

\def\S{\mathbb S}

\title{Local chromatic number of quadrangulations of surfaces}

\author{
Bojan Mohar\thanks{Supported in part by an NSERC Discovery Grant,
  by the Canada Research Chair program, and by Research Grant P1--0297 of ARRS, 
  Slovenia.}~\thanks{On leave from:
  IMFM \& FMF, Department of Mathematics, University of Ljubljana, Ljubljana,
  Slovenia.}\\
  Department of Mathematics\\ Simon Fraser University\\ 
  Burnaby, B.C.\\ \texttt{mohar@sfu.ca} 
\and
G\'abor Simonyi\thanks{Research partially supported by the
  Hungarian Foundation for Scientific Research Grant (OTKA) Nos.\ K76088 and
  NK78439.}\\
  Alfr\'ed R\'enyi Institute of Mathematics\\
  Hungarian Academy of Sciences\\ Budapest, Hungary\\
  \texttt{simonyi@renyi.hu}
\and
G\'abor Tardos\thanks{Supported in part by an NSERC Discovery Grant, 
  by the Canada Research Chair program, 
  and by Hungarian Foundation for Scientific Research Grant (OTKA) 
  no.\ 78439.}\\
  Alfr\'ed R\'enyi Institute of Mathematics\\
  Hungarian Academy of Sciences, Budapest, Hungary\\ 
  and \\ Department of Computer Science\\ 
  Simon Fraser University, Burnaby, B.C.\\ \texttt{tardos@cs.sfu.ca}
}

\date{}
\begin{document}
\maketitle

\begin{abstract}
The {\em local chromatic number} of a graph was introduced in \cite{EFHKRS}.
In \cite{ST} a connection to
topological properties of (a box complex of) the graph was established and in
\cite{STV} it was shown that if a graph is {\em strongly topologically
$4$-chromatic} then its local chromatic number is at least four. As a
consequence one obtains a generalization of the following theorem of Youngs 
\cite{Y}: If a quadrangulation of the projective plane is not bipartite it has
chromatic number four. The generalization states that in this case the local
chromatic number is also four.

Both papers \cite{AHNNO} and \cite{MS} generalize Youngs's result to arbitrary
non-orientable surfaces replacing the condition of the graph being not
bipartite by a more technical condition of an {\em odd} quadrangulation. This
paper investigates when these general results are true for the local chromatic
number instead of the chromatic number. Surprisingly, we find out that (unlike
in the case of the chromatic number) this depends on the genus of 
the surface. For the non-orientable surfaces of genus at most four, 
the local chromatic number of any odd quadrangulation is at least four,
but this is not true for non-orientable surfaces of genus 5 or higher.

We also prove that face subdivisions of odd quadrangulations and 
Fisk triangulations of arbitrary surfaces exhibit the same behavior for the
local chromatic number as they do for the usual chromatic number.
\end{abstract}

\section{Introduction}

We start by defining the terms in the title. By a {\em surface} we mean a
compact connected 2-manifold without boundary. A {\em quadrangulation} of a
surface is a loopless graph on a surface with all faces being quadrilaterals. 
We allow parallel edges and quadrilateral faces with some of its vertices or
even edges coinciding. Given an arbitrary
orientation of all faces we can count the edges that break consistency of
the orientation of the surface. (We say that consistency is not broken at some
edge if, considering the orientation of a face as a closed walk along its 
boundary, we traverse the edge in opposite directions when considering 
the faces at its two sides.)  
Note that reversing the orientation of a face changes the status
of its four edges (except for the edges appearing twice on the boundary)
and thus the parity of the number of edges breaking consistency does not
change. We call a quadrangulation {\em even} or {\em odd\/} depending on this
parity. 
By the previous remark, the parity of the quadrangulation is determined by the
quadrangulation itself and is independent of the actual orientation of the
faces.  

Note that every quadrangulation
of an orientable surface is even as all faces can be oriented consistently.

We denote by $\chi(G)$ the {\em chromatic number} of the graph $G$.
The following result was proved independently by Archdeacon et al.\
\cite{AHNNO} and Mohar and Seymour \cite{MS}. 

\begin{thm}[\cite{AHNNO,MS}]\label{am}
For an odd quadrangulation $G$ of a surface we have $\chi(G)\ge4$.
\end{thm}

The theorem only applies to non-orientable surfaces as orientable ones
have no odd quadrangulations. Youngs \cite{Y} has established this earlier for
the case of the projective plane. This result is simpler to state, since 
being even and being bipartite is equivalent for quadrangulations 
of the projective plane and since
every quadrangulation of the projective plane is 4-colorable.

\begin{thm}[\cite{Y}]The chromatic number of a quadrangulation of the
projective plane is either two or four.
\end{thm}

In this paper we generalize the above results for the local chromatic number.
The local chromatic number of a graph is defined in \cite{EFHKRS} as the
minimum number of colors that must appear within distance $1$ of a vertex. 
For the formal definition, let $N(v)=N_G(v)$ denote the {\em neighborhood} of a
vertex $v$ in a graph $G$, that is, $N(v)$ is the set of vertices 
adjacent to $v$. 
For an integer $k\ge1$, we call a proper vertex-coloring $c$ of a graph $G$ a
{\em local $k$-coloring} if $|c(N(v))|\le k-1$ for every vertex $v$ of $G$.
The {\em local chromatic number} $\psi(G)$ of $G$ is the smallest $k$ such
that $G$ has a local $k$-coloring.

The $-1$ term comes traditionally from considering ``closed neighborhoods''
$N(v)\cup\{v\}$ and results in a simpler form of the relations with other
coloring parameters, like the trivial bound $\psi(G)\le\chi(G)$.

More on the local chromatic number, including examples of graphs with
arbitrarily high chromatic number and local chromatic number $3$, can be
found in \cite{EFHKRS}, see also \cite{ST, STV}.

First we state the generalization of Youngs's theorem. This result was
hinted in Remark~3 of \cite{STV}, here we give a full proof for the sake of
completeness in the next section. 

\begin{thm}\label{pplane}
The local chromatic number of a quadrangulation of the
projective plane is two or four.
\end{thm}

Note that the non-orientable surfaces are determined by a positive invariant,
their {\em genus}. The projective plane has genus one. Next we generalize
Theorem~\ref{pplane} for higher genus surfaces. This is the main
result of the present paper. Surprisingly, and unlike in
the case of the chromatic number, the situation depends on the genus of the
surface.

\begin{thm}\label{main}
\begin{description}
\item[(i)] If $G$ is an odd quadrangulation of a (non-orientable) surface of
genus at most four, then $\psi(G)\ge4$.
\item[(ii)] Every non-orientable surface of genus at least five admits an odd
quadrangulation that has a local $3$-coloring using six colors.
\end{description}
\end{thm}

We remark that at least five colors are needed to locally $3$-color odd 
quadrangulations or in fact any graph with chromatic number higher than
$3$. This is a consequence of the fact that any local $3$-coloring  
with four colors can easily be transformed into a $3$-coloring by changing the
color of each vertex colored 4 to the color 1, 2, or 3, which is not 
used in its neighborhood.

If we insist on five colors, the threshold lies higher:

\begin{thm}\label{5colors}
\begin{description}
\item[(i)]
If $G$ is an odd quadrangulation of a (non-orientable) surface of genus at
most six, then $G$ has no local $3$-coloring with at most five colors.
\item[(ii)]
Every non-orientable surface of genus at least seven admits an odd
quadrangulation that has a local $3$-coloring using five colors.
\end{description}
\end{thm}

For (some special cases of) Theorem~\ref{main}(i) we give several proofs
because they use very different approaches.
In Section~\ref{sec:box} we present the argument hinted in \cite{STV}
for Theorem~\ref{pplane}, i.e., for the projective plane. It is based on
the concept of the hom space of a graph, which is just a slightly different
version of the graph complexes known as box complexes, see \cite{MZ}, or those
more generally called hom complexes, see (the cited second edition of)
\cite{Mat}.  
In Section~\ref{sec:alg} we give an algebraic proof
for the non-existence of odd locally $3$-chromatic quadrangulations that works
for both the projective plane and the Klein bottle. In
Section~\ref{sec:comb} we give a combinatorial argument that works for all
surfaces of genus at most four, and henceforth proves Theorem~\ref{main}(i).
The results provided in Sections~\ref{sec:box} and \ref{sec:alg} follow
also from the result of Section~\ref{sec:comb}. But both of the former sections
indicate connections between different concepts that are not shown by the
combinatorial argument in Section~\ref{sec:comb}. In Section~\ref{sec:box} we
use a general connection between any surface and the hom space of its
quadrangulation. In Section~\ref{sec:alg} we show how an algebraic question
about so-called semi-free groups is related to our graph theoretic
problems. As a consequence of this connection and the existence of the
quadrangulations claimed in Theorem~\ref{main} we can partially answer the
algebraic question.

The locally $3$-chromatic quadrangulations claimed in Theorems~\ref{main}(ii)
and \ref{5colors}(ii) are
constructed in Section~\ref{sec:const}. These constructions are motivated by
the argument in the preceding section. In Section~\ref{sec:5colors} we prove
Theorem~\ref{5colors}(i).

By inserting a vertex in each face (and joining it to the vertices on the
boundary of that face) of an odd quadrangulation we obtain a triangulation
that exhibits similar behavior for the chromatic number: any subgraph
contained in a contractible part of the surface is 3-colorable, but the
chromatic number of the whole graph is at least five. It is shown in
Section \ref{sect:triang} that the same holds for the local chromatic number.
When coloring graphs on surfaces, another family of triangulations exhibits
unusual behavior \cite{MT}; these are triangulations in which all vertices
except two have even degree, and the two vertices of odd degree are adjacent.
Fisk \cite{Fi} proved that such triangulations cannot be 4-colored, and we show
that the same holds for the local chromatic number. 
Finally, Section~\ref{sec:last} contains related remarks and observations.

\section{Hom spaces}\label{sec:box}

In this section we prove Theorem~\ref{pplane} for which only a very rough
argument was presented in \cite{STV}. We start with the definition of medial
graphs that will be used in the sequel.

Let $G$ be a graph on a surface $S$. We define the {\em medial graph} $M(G)$,
also embedded in $S$, as follows. The vertices of $M(G)$ correspond to
edges of 
$G$. For an edge $e$ of $G$ we choose an interior point $v_e$ of $e$ as the
corresponding vertex of $M(G)$. For each vertex $x$ of $G$ we consider the edges
of $G$ incident to $x$ and place a cycle on the corresponding vertices of
$M(G)$ as follows. 
Let $e_1,\ldots,e_d$ be the edges incident to $x$ listed in the cyclic
order these edges leave $x$. For $1\le i<d$ we connect $v_{e_i}$ and
$v_{e_{i+1}}$ with an edge of $M(G)$ drawn inside the face of $G$ bounded by
$e_i$ and $e_{i+1}$. We connect $v_{e_d}$ and $v_{e_1}$ the same way.

Notice that $M(G)$ is a 4-regular graph embedded in $S$ that
has two types of faces. The {\em star faces} are the ones
containing a vertex $x$ of $G$ and bounded by the cycle we introduced on the
vertices of $M(G)$ corresponding to the edges of $G$ incident to $x$ (the star
of $x$). The {\em cycle faces} are the remaining faces of $M(G)$: each is
contained in a single face of $G$ and its vertices correspond to the edges
on the corresponding facial walk.

The proof of Theorem~\ref{pplane} is based on the notion of the {\em hom
space} $H(G)$ defined for any finite graph $G$. For the precise definition of
this and related concepts see, e.g., \cite{ST}. Here we note that $H(G)$
is a 
cell complex, whose cells correspond to complete bipartite subgraphs of $G$.
More precisely, the cells are of the form $A\uplus B$, where $A$ and $B$ are
non-empty, disjoint sets of vertices, such that every vertex of $A$ is connected in
$G$ to every vertex of $B$. (We use the notation $A\uplus B$ for the ordered
pair containing $A$ and $B$ as it is customary in this context,
cf.\ \cite{Mat}.) The vertices of $H(G)$ are of the form
$\{x\}\uplus\{y\}$, where $xy$ is an edge of $G$. 
In particular, every edge $xy$ gives rise to two vertices, $\{x\}\uplus\{y\}$
and $\{y\}\uplus\{x\}$, in $H(G)$. The map $\sigma$ that
switches $A\uplus B$ with $B\uplus A$ is a continuous involution that makes
$H(G)$ into a $\2$-space. (For the notion and basic properties of $\2$-spaces
we also refer to \cite{Mat}.)

Next we show an interesting topological connection between a surface and the
hom space of its quadrangulation from which the proof of
Theorem~\ref{pplane} follows easily. 

Let $G$ be a quadrangulation of a surface $S$. Let $\alpha$  be a map that
maps the vertex 
$v_e$ of $M(G)$ to the class of $\{x\}\uplus\{y\}$ in $H(G)/\sigma$ where
$e=xy$. Now $\alpha$ naturally extends to the edges of $M(G)$ (mapping
them to 1-cells connecting the corresponding vertices) and also to the faces
of $M(G)$ as follows. The image of the cycle face of $M(G)$ 
corresponding to the face $xyzt$ of
$G$ will be the class of the cell $\{x,z\}\uplus\{y,t\}$. The image of the
star face containing the vertex $v$ of $G$ is inside the class of the cell
$\{v\}\uplus N(v)$. This makes $\alpha$ a continuous map $\alpha:S\to
H(G)/\sigma$.

The map $\alpha$ lifts to a map $\beta:T\to H(G)$ where $T$ is 
a double cover of $S$. 
(See \cite{Mas} for more details about covering spaces of graphs and surfaces.)
Let us observe that the double cover restricted to the graph $G$ of the
quadrangulation is isomorphic to the categorical product $G\times K_2$
(sometimes 
called the Kronecker product or the direct product), and it is connected if
and only if $G$ is not 
bipartite. Let us consider the involution on $T$ interchanging the pairs of points 
with the same image in $S$. This makes $T$ into a $\2$-space and $\beta$ a
$\2$-map. Notice that $T$ is the union of two disjoint copies of $S$ if and
only if $G$ is bipartite.

\medskip

{\em Remark 1.}
Notice that the parity of the length of a cycle of a quadrangulation of a
surface is determined by the homotopy type of the cycle as a surface
cycle. These parities for the different homotopy types determine the double
cover $T$ as follows: an even cycle on $S$ lifts to two cycles on $T$, while
an odd cycle on $S$ lifts to a single cycle on $T$.

If we assume that all $4$-cycles in $G$ are faces of the quadrangulation and 
$G$ is not a complete bipartite graph $K_{2,i}$ ($i\le3$), then
$H(G)$ and $T$ are homotopy equivalent, with $\beta$ providing one direction of
the homotopy equivalence. The exceptional graphs $K_{2,i}$ for $i\ge1$ 
quadrangulate the sphere and for $i\le3$ all $4$-cycles are faces. 
To see the homotopy equivalence consider the maximal cells of $H(G)$---these
correspond to maximal complete bipartite subgraphs of $G$, which are the
face cycles and the stars of vertices. It's not hard to see that $G$ has no
vertex of degree less than $3$. If $G$ is
$3$-regular, $T$ and $H(G)$ are actually homeomorphic and $\beta$ can be
chosen to be a homeomorphism. Two $(d-1)$-dimensional simplicial cells of
$H(G)$ correspond to the star of a degree $d$ vertex of $G$, so if $G$ is not
$3$-regular, $T$ and $H(G)$ are not homeomorphic. It is easy to see though that
only the edges of a $d$-cycle on the boundary of these cells appear in other
cells, therefore the higher dimensional cell can be collapsed to a $d$-gon
without changing the homotopy type. After this collapse of every higher
dimensional cell we obtain a complex homeomorphic to $T$. 
\hfill$\Diamond$

\bigskip

\noindent{\bf Proof of Theorem~\ref{pplane}:}
The upper bound $4$ follows from $\psi(G)\le\chi(G)$ and from the easy part of
the Youngs theorem (that every quadrangulation of the projective plane
is 4-colorable). The nontrivial part of the result
is establishing that if $G$ is not bipartite it must have no local 3-coloring. 

Consider the map $\beta:T\to H(G)$ constructed above. Here $T$ is the double
cover of the projective plane and it is not the union of two disjoint copies
of that plane as $G$ is not bipartite. The only possibility left for $T$ is
therefore the sphere $T=\S^2$. The existence of a $\2$-map $\beta:\S^2\to H(G)$
is signified by saying $G$ is strongly topologically $4$-chromatic and for
such graphs it is established in \cite{STV} that $\psi(G)\ge4$.\qed

\medskip

{\em Remark 2}: Among the graphs investigated in \cite{ST, STV} $4$-chromatic
generalized Mycielski graphs are known to quadrangulate the projective plane,
see \cite{GyJS}, while deleting some edges from $4$-chromatic Schrijver graphs
one obtains quadrangulations of the Klein bottle. In fact, it was the
investigation of the local chromatic number of the latter graphs that led us
to the proof presented in the next section and to consider the local chromatic
number of surface quadrangulations in general. \hfill$\Diamond$

\section{Semi-free groups}\label{sec:alg}

In this section we present an intermediate step toward Theorem~\ref{main}(i).
This algebraic approach proves
the statement for the non-orientable surfaces of genus one or two, i.e., for
the projective plane and the Klein bottle.

Let $H$ be a graph. The {\em semi-free} group $\Gamma_H$ is the following
group given by generators and relations. The generators are the vertices of $H$
while the relations are $xy=yx$ for the edges $xy$ of $H$. If $E(H)=\emptyset$
we obtain the free group, for a complete graph we obtain a free Abelian group. 
Thus, one can consider semi-free groups as a common generalization of
free and free Abelian groups.

Let $G$ be a quadrangulation of a surface $S$ and let $c$ be a local
$3$-coloring of $G$. Our goal is to prove a lower bound on the genus of
$S$ if the quadrangulation is odd.

Let $m$ be the number of colors used by the coloring $c$ and assume that these
are the elements of $[m]=\{1,\ldots,m\}$. Consider the
Kneser graph $H={\mathrm KG}(m,2)$ and the semi-free group
$\Gamma=\Gamma_H$. The generators of $\Gamma$ are the vertices of $H$, which
are the $2$-element subsets of $[m]$. (Recall that two such vertices are
adjacent in ${\mathrm KG}(m,2)$ iff they represent disjoint subsets.) For
$i,j\in[m]$ we introduce the notation 
$$
    x_{i,j}=\left\{\begin{array}{lll}
               \ve&\;&\hbox{if }i=j\\
               \{i,j\}&&\hbox{if }i<j\\
               \{i,j\}^{-1}&&\hbox{if }j<i.
            \end{array}\right.
$$
where $\ve$ is the identity element of $\Gamma$ and $\{i,j\}$ is a generator of
$\Gamma$ if $i\ne j$. These group elements clearly satisfy
$x_{i,j}=x_{j,i}^{-1}$ for any $i,j\in[m]$ and by the definition of the Kneser
graph we have $x_{i,j}x_{k,l}=x_{k,l}x_{i,j}$ whenever $\{i,j\}$ and $\{k,l\}$
are disjoint. Moreover, we have
\begin{equation}\label{eq}
x_{i,j}x_{j,k}=x_{i,k}\;\hbox{ whenever }|\{i,j,k\}|\le2.
\end{equation}

Consider the medial graph $M(G)$ of $G$. (See the definition in
Section~\ref{sec:box}.) We consider the edges of $M(G)$ as oriented edges with
both orientations of each (unoriented) edge being present. We label these
oriented edges by elements of the group $\Gamma$. Consider the oriented edge
$w=(v_e,v_f)$, where $e$ and $f$ are edges of $G$. Notice that $e=ab$ and 
$f=ad$ must be adjacent edges. In this case we label the oriented edge
$w$ by $l(w)=x_{c(b),c(d)}\in\Gamma$. We let the label $l(W)$ of a
walk $W$ on $M(G)$ be the product of the labels along the walk.

\begin{figure}[htb]
\epsfxsize6cm \centerline{\epsffile{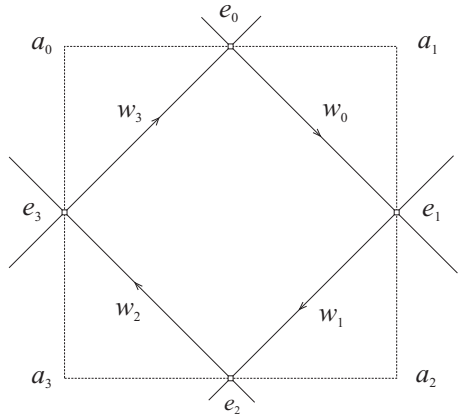}} 
\caption{A cycle face in $M(G)$}
\label{fig:1}
\end{figure}

\begin{lem}\label{label}
For the opposite orientations $w$ and $w'$ of the same edge we
have $l(w')=l(w)^{-1}$. For a walk\/ $W$ around a face of\/ $M(G)$ we have
$l(W)=\ve$.
\end{lem}

\proof The first statement is trivial.
A walk $W$ around a cycle face of $M(G)$ consists of four edges. See
Figure~\ref{fig:1} for this case. Let the vertices around the face of $G$ 
containing $W$ be $a_0,a_1,a_2,a_3$, connected by the edges $e_i=a_ia_{i+1}$
(where all indices are taken modulo 4).
The walk $W$ consists of the oriented edges $w_i=(e_i,e_{i+1})$ \ ($0\le
i\le3$) and we have $l(w_i)=x_{c(a_i),c(a_{i+2})}$. Here
$l(w_1)=x_{c(a_1),c(a_3)}=x_{c(a_3),c(a_1)}^{-1}=l(w_3)^{-1}$. Similarly,
we have $l(w_0)=l(w_2)^{-1}$. The labels $l(w_0)=x_{c(a_0),c(a_2)}$ and
$l(w_1)=x_{c(a_1),c(a_3)}$ commute since if the sets $\{c(a_0),c(a_2)\}$ and
$\{c(a_1),c(a_3)\}$ were not disjoint, $c$ would not be a proper coloring.
Consequently, $l(W)=\ve$.

\begin{figure}[htb]
\epsfxsize7cm \centerline{\epsffile{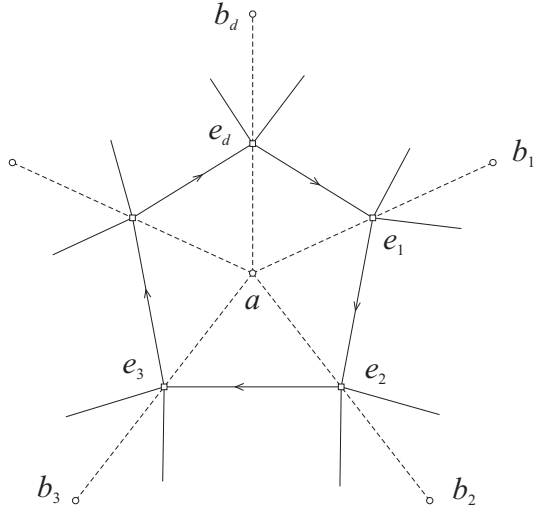}} 
\caption{A star face in $M(G)$}
\label{fig:2}
\end{figure}

Let us now consider the walk $W$ around a star face of $M(G)$ as on 
Figure~\ref{fig:2}. Let $a$
be the vertex of $G$ inside this face and $b_1,\ldots,b_d$ the neighbors of
$a$ in this order. We have
$$l(W)=x_{c(b_1),c(b_2)}x_{c(b_2),c(b_3)}\ldots
  x_{c(b_{d-1}),c(b_d)}x_{c(b_d),c(b_1)}.$$
As $c$ is a local $3$-coloring there are at most $2$ distinct ones among 
the colors $c(b_i)$. So Equation~(\ref{eq}) applies and simplifies the above 
expression to $l(W)=x_{c(b_1),c(b_1)}=\ve$.\qed

\bigskip

Let us fix a vertex $v_{e_0}$ of $M(G)$ as the base point of $S$. If two
closed walks on $M(G)$ starting at the base point are homotopic on $S$, then
one can be transformed to the other by a sequence of steps, where each step is
the introduction or removal of a closed sub-walk which is either a facial walk
or an edge traversed both ways. By our last lemma these steps do not alter the
label of the walk, so the label is determined by the homotopy type. As any
homotopy type can be represented by a closed walk on $M(G)$, we have a map
$\alpha:\pi_1(S)\to\Gamma$ from the fundamental group $\pi_1(S)$ of $S$ 
that gives as its value to each element of the fundamental group 
the label of the walks in the corresponding homotopy class. Clearly, this is a
homomorphism. 

We shall need another property of products in the semi-free group
$\Gamma_{{\mathrm KG}(m,2)}$.

\begin{lem}\label{lem:new}
Let $c_1,\dots,c_t$ be colors in $[m]$, where $c_{i+1}\ne c_i$ for
$i=1,\dots,t-1$, and $c_1\ne c_t$. 
Let $\alpha=x_{c_t,c_2}x_{c_1,c_3}\ldots x_{c_{t-2},c_t}x_{c_{t-1},c_1}$.
If\/ $t$ is odd, then $\alpha\ne\ve$.
\end{lem}

\proof
Let us assume the lemma is false and consider a counterexample where a minimum
number among the factors $X_i=x_{c_{i-1},c_{i+1}}$ are non-identity. 
(Here and in the rest of the proof, all indices are understood modulo $t$.)
There must be such non-identity factors as otherwise we would have
$c_{i-1}=c_{i+1}$ for all $i$ and therefore (as $t$ is odd) that all the $c_i$
are the same.

Note that each factor $X_i$ is the identity or a generator or the inverse of a
generator in the semi-free group $\Gamma$. We use the observation of Baudisch
\cite{Baud} that if a nontrivial product of generators and their inverses in a
semi-free group is the identity, then one finds a generator and its inverse
in the product such that they commute with all factors separating them (so they
cancel each other). We may therefore choose indices $1\le i<j\le t$ with
$X_i\ne \ve$, $X_j=X_i^{-1}$ and $X_i$ commuting with all $X_s$ for
$i<s<j$. We choose $j$ so that $X_j$ is the first occurrence of $X_i^{-1}$
after the factor $X_i$. We have $c_s=c_{i+1}$ for
all $i<s<j$ with $s-i$ odd, because otherwise for the smallest $s$ breaking
this rule either $X_{s-1}$ would not commute with $X_i$ or $X_{s-1}=X_i^{-1}$
would hold contradicting the minimality of $j$. For $i<s<j$ with $s-i$ even we
have $c_s\ne c_{i+1}$, because adjacent colors must differ. This implies that
$j-i$ is even as $c_{j-1}=c_{i+1}$. Furthermore, $c_s\ne c_{i-1}$ for $i<s<j$
as otherwise $X_{s-1}$ and $X_i$ would not commute for the first such $s$.

Let us modify now the colors $c_s$ for $i<s<j$, $s-i$ odd to
$c'_s=c_{i-1}$ 
without changing the other colors. This does not create equal neighboring
colors and changes only two factors: $X_i'=X_j'=\ve$. We obtain another
counterexample to the claim this way, and this counterexample has a smaller
number of non-identity factors. This contradicts our choice and by this
contradiction proves the lemma.\qed

\begin{lem}\label{odd}
If $\alpha(y)=\ve$ for some class $y$ in the fundamental group of $S$, then any
walk on $G$ representing $y$ has even length.
\end{lem}

\proof Assume for contradiction that $\alpha(y)=\ve$ but the length $t$ of 
a closed walk $W$ representing $y$ is odd. 
Note that here we have the non-standard view of the walk $W$ as 
``starting'' at the base point $v_{e_0}$ that happens to be in the
middle of the edge $e_0$ of $G$. However, this means that the walk $W$
uses at the beginning one half of the edge $e_0$, and traverses the other 
half when coming back at the end. To apply
the definition of $\alpha$, we have to find first a walk $W'$ on $M(G)$
homotopic to $W$. Let $a_1,\ldots,a_t$ be the vertices of $G$ along
$W$. Since $W$ starts and ends at the base point, we have
$e_0=a_1 a_t$. Let $e_i$ be the edge of $G$ connecting $a_i$ and $a_{i+1}$
for $1\le i<t$. We construct the walk $W'$ on $M(G)$ as follows. It starts at
the base point $v_{e_0}$ and then it passes through all the
points $v_{e_i}$ for $1\le i<t$ in order before returning to the base
point. We let the part $W_i$ of $W'$ between $v_{e_{i-1}}$ and $v_{e_i}$
follow the boundary of the star face of $M(G)$ containing $a_i$ for $1\le i\le
t$. The indices here and in the rest of this proof are understood modulo $t$.

Let $c_i=c(a_i)$. As in the proof of Lemma~\ref{label}
we can use the fact that $c$ is a local 3-coloring and Equation~(\ref{eq}) to
conclude that $l(W_i)=x_{c_{i-1},c_{i+1}}$. We thus have
$$\alpha(y)=l(W')=\prod_{i=1}^tl(W_i)=x_{c_t,c_2}x_{c_1,c_3}\ldots
x_{c_{t-2},c_t}x_{c_{t-1},c_1}.$$
Lemma \ref{lem:new} shows that a product of this form with
$t$ odd is never the identity if $c_i\ne c_{i+1}$ for all $1\le i\le
t$. Note that neighboring colors $c_i$ and $c_{i+1}$ are distinct in our case
as $c$ is a proper coloring, so this completes the proof.
\qed

\bigskip
{\em Remark 3.} One can try to simplify the argument presented in
this section by letting all pairs of generators commute in $\Gamma$, i.e.,
considering a commutative factor of $\Gamma$. The argument breaks down because
Lemma~\ref{odd} does not hold in this case. Assume $G$ has a cycle of length
$9$ with the vertices along the cycle colored $2,1,2,3,1,3,4,1,4$.
Lemma~\ref{odd} applies and $\alpha(y)\ne \ve$ in $\Gamma$ for the class $y$
containing this cycle. But $\alpha(y)$ vanishes in any commutative factor of
$\Gamma$.
\hfill$\Diamond$

\begin{thm}\label{semifree}
Let $g$ be a positive integer and assume the equality $z_1^2z_2^2\ldots
z_g^2=\ve$ implies $z_1z_2\ldots z_g=\ve$ in the semi-free group
$\Gamma_{{\mathrm KG}(m,2)}$, where, as before, $\ve$ denotes the identity
element of the group.  
Then no odd quadrangulation of a non-orientable surface of genus $g$ has a local 
$3$-coloring using at most $m$ colors.
\end{thm}

\proof 
Let $S$ be the non-orientable surface of genus $g$. The fundamental group
$\pi_1(S)$ can be presented as $<y_1,\ldots,y_g\mid y_1^2y_2^2\ldots y_g^2=1>$. 
(Cf., e.g., \cite{Mas} for more details about fundamental groups of surfaces.)
Assume $S$ has a quadrangulation $G$ with a local $3$-coloring using $m$ colors. 
We need to prove that this quadrangulation is even. With
the homomorphism $\alpha:\pi_1(S)\to\Gamma_{{\mathrm KG}(m,2)}$ defined above we
let $z_i=\alpha(y_i)$ for $1\le i\le g$. We clearly have $z_1^2z_2^2\ldots
z_g^2=\ve$. By our assumption this implies $z_1z_2\ldots z_g=\ve$, so
$\alpha(y)=\ve$ for $y=y_1y_2\ldots y_g$.

The surface $S$ can be obtained by placing $g$ crosscaps on the sphere and
$y_i$ can be chosen to be the class of a loop going only through 
the $i$th crosscap, so $y$ will be the class of a loop going
through each crosscap once.

If we cut $S$ along a simple cycle in the homotopy class $y$ we obtain an
orientable surface. If $y$ can be represented as a cycle $C$ of $G$, then one
can cut $S$ along $C$ and consistently orient the
obtained surface and in it all the faces of the quadrangulation $G$. With this
orientation exactly the edges of $C$ break the consistency of the orientation
and thus the parity of the quadrangulation $G$ is the same as the parity
of the length of the cycle $C$.

In general we cannot assume the existence of a simple cycle in $G$
representing the class $y$ but we can always represent $y$ by a closed walk
$W$ on $G$. We claim that the
length of the walk has the same parity as the quadrangulation $G$. To see
this one can refine the quadrangulation without changing its parity or that of
any cycle till one finds a simple cycle homotopic to $W$ or alternatively one
can prove that there is an orientation of the faces of $G$ where the
consistency is broken at exactly the edges traversed an odd number of times by
$W$.

By Lemma~\ref{odd}, the length of $W$ must be even because $\alpha(y)=\ve$. This
shows that $G$ is an even quadrangulation and completes the proof of the
theorem.\qed

\bigskip

We proved Theorem~\ref{semifree} to find another proof for
Theorem~\ref{pplane} that extends also to higher
genus surfaces.

The non-existence of an odd, locally $3$-colorable
quadrangulation of the projective plane follows now from the fact that $x^2=\ve$
implies $x=\ve$ in every semi-free group. This statement, and more generally
that semi-free groups have no torsion elements, was proved by Baudisch
\cite{Baud}.

In another paper \cite{Baud2}, Baudisch proved that any two non-commuting
elements of a semi-free group freely generate a free group of rank 2. 
Thus $x^2y^2=\ve$
implies that $x$ and $y$ commute. So we have $(xy)^2=x^2y^2=\ve$ and by the
earlier result $xy=\ve$. This establishes that the Klein bottle (the 
non-orientable surface of genus 2) has no odd, locally $3$-colorable
quadrangulations.

To apply Theorem~\ref{semifree} to the next surface we would need that
$x^2y^2z^2=\ve$ implies $xyz=\ve$ in semi-free groups, but we were not able to
prove this.

We can turn Theorem~\ref{semifree} around and use its counterpositive form. 
From the odd quadrangulation claimed in Theorem~\ref{main}(ii) we conclude
that 
there are elements $z_1$, $z_2$, $z_3$, $z_4$ and $z_5$ in the semi-free group
$\Gamma_{{\mathrm KG}(6,2)}$ satisfying $z_1^2z_2^2z_3^2z_4^2z_5^2=\ve$ and 
$z_1z_2z_3z_4z_5\ne \ve$. With ``reverse
engineering'' the construction in Section~\ref{sec:const}, one can actually
find these elements, see Table~1. Note that as a result of the properties of
the construction our example uses only nine of the fifteen generators of
$\Gamma_{{\mathrm KG}(6,2)}$, so our example actually lives in the semi-free
group $\Gamma_{C_3^2}$, where $C_3^2$ is the square of the 3-cycle.
\medskip

\begin{eqnarray*}
z_1&=&a_{25}a_{14}a_{24}^{-1}a_{15}^{-1}\\
z_2&=&a_{15}a_{24}a_{14}^{-1}a_{16}a_{36}^{-1}a_{15}^{-1}a_{35}a_{14}a_{24}^{-1}a_{15}^{-1}\\
z_3&=&a_{15}a_{24}a_{35}^{-1}a_{36}a_{14}^{-1}a_{26}^{-1}\\
z_4&=&a_{26}a_{14}a_{36}^{-1}a_{24}^{-1}a_{34}a_{14}^{-1}\\
z_5&=&a_{14}a_{34}^{-1}a_{36}a_{24}a_{16}^{-1}a_{15}a_{25}^{-1}a_{14}^{-1}\\
z_1^2z_2^2z_3^2z_4^2z_5^2&=&\ve\\
z_1z_2z_3z_4z_5&=&a_{25}a_{16}a_{15}^{-1}a_{14}a_{16}^{-1}a_{15}a_{25}^{-1}a_{14}^{-1}\ne\ve
\end{eqnarray*}
\centerline{{\bf Table 1.} Five elements in the semi-free group
$\Gamma_{{\mathrm KG}(6,2)}$ showing peculiar behavior.}
\centerline{The generator corresponding to the vertex $\{i,j\}$ of ${\mathrm
KG}(6,2)$ is denoted by $a_{ij}$.}
\medskip

From the quadrangulation, whose existence is claimed in
Theorem~\ref{5colors}(ii), one can construct a similar list of $7$
elements in the semi-free group corresponding to ${\mathrm KG}(5,2)$, that is, 
the Petersen graph. See Table~2.
\medskip

\begin{eqnarray*}
w_1&=&b_{23}b_{13}^{-1}b_{24}^{-1}b_{14}\\
w_2&=&b_{14}^{-1}b_{13}b_{24}b_{35}b_{25}^{-1}\\
w_3&=&b_{25}b_{24}^{-1}b_{35}^{-1}b_{34}\\
w_4&=&b_{34}^{-1}b_{35}b_{24}b_{34}^{-1}b_{13}^{-1}b_{35}^{-1}b_{15}b_{23}^{-1}b_{35}^{-1}b_{34}\\
w_5&=&b_{34}^{-1}b_{35}b_{15}^{-1}b_{23}b_{12}b_{34}\\
w_6&=&b_{34}^{-1}b_{12}^{-1}b_{23}^{-1}b_{15}b_{45}^{-1}\\
w_7&=&b_{23}b_{45}a_{15}^{-1}b_{34}b_{13}b_{23}^{-1}\\
w_1^2w_2^2w_3^2w_4^2w_5^2w_6^2w_7^2&=&\ve\\
w_1w_2w_3w_4w_5w_6w_7&=&b_{23}b_{35}b_{34}^{-1}b_{13}^{-1}b_{35}^{-1}b_{34}b_{13}b_{23}^{-1}\ne\ve
\end{eqnarray*}
\centerline{{\bf Table 2.} Seven elements in the semi-free group
$\Gamma_{{\mathrm KG}(5,2)}$ showing peculiar behavior.}
\centerline{The generator corresponding to the vertex $\{i,j\}$ of ${\mathrm
KG}(5,2)$ is denoted by $b_{ij}$.}
\medskip

\section{A combinatorial approach}
\label{sec:comb}

In this section we present the most elementary and, at the same time, at least
by our current knowledge, the most
effective approach toward proving that certain quadrangulations have no local
3-colorings. It is based on examining properties of minimal counterexamples.

Let us call a triple $(G,S,c)$ a {\em suitable quadrangulation}, if $G$ is an
odd quadrangulation of the (non-orientable) surface $S$ and $c$ is a local 
$3$-coloring of $G$. A suitable quadrangulation $(G,S,c)$ is called a 
{\em minimal quadrangulation} if no
surface of genus less than that of $S$ has a suitable quadrangulation and $S$
has no suitable quadrangulation with fewer faces.

\begin{lem}\label{reduction}
If $(G,S,c)$ is a minimal quadrangulation and $F$ is a face of this
quadrangulation, then $F$ has four distinct vertices. Moreover, these vertices
receive two or four distinct colors in the local\/ $3$-coloring. 
The neighbors of any vertex in $G$ receive exactly two distinct colors.
\end{lem}

\proof As $c$ is a proper coloring, we have no loops, so only diagonally
opposite vertices of $F$ can coincide. Let us assume that the same vertex $x$
appears at both endpoints of the diagonal $d$ of the face $F$. 
Then the two edges of $F$ connecting $x$ to another vertex of $F$ are parallel
edges or they coincide. Let
us cut the face $F$ out from $S$ and close $S$ up by identifying these pairs
of parallel edges of $F$ or by removing the edge if the two
neighboring sides of $F$ coincided. Let $S'$ be the space obtained from
$S$ this way and $G'$ be the resulting graph on $S'$. As we only identified
parallel edges, $c$ is a local $3$-coloring of $G'$. Considering
any orientation of the faces of $G$ and the same orientation of the
faces of $G'$ one sees that an edge of G' obtained by identification breaks
consistency of this orientation if and only if exactly one of the
corresponding two edges in $G$ does so, other edges of $G'$ break
consistency of the orientation in $G'$ if and
only if they do so in $G$, while edges removed from $G$ did not break
consistency there. Therefore, since $G$ is an odd quadrangulation,
$G'$ must be odd, too. Note, 
however, that $S'$ is not necessarily a surface as the neighborhood of $x$ can
be strange (see below), so for the
previous sentence to make sense we have to allow a somewhat 
extended definition of parity of quadrangulations.

In the following case analysis we always find a quadrangulation showing that
$(G,S,c)$ is not minimal. This is done in slightly different ways depending on
the topology of the diagonal $d$.

In the simplest case $d$ is a one-sided simple closed curve on 
$S$ (which is clearly non-separating). In this case $S'$ is a surface. (To see
this, it is enough to check how the faces incident to $x$ are
arranged around this point.) So
$(G',S',c)$ is a suitable quadrangulation of $S'$. 
The surface $S'$ is non-orientable as it has an odd
quadrangulation. Its Euler characteristic is one more than that of $S$, as $G'$
has one fewer face and two fewer edges, and it has the same number of
vertices as $G$. So the genus of $S'$ is one less than that of $S$
contradicting the minimality of $(G,S,c)$.

Our second case is when a pair of coinciding neighboring edges got removed
from $G$. This makes $d$ a (trivial) separating cycle on $S$, $S'$ a
surface homeomorphic to $S$ and $(S',G',c)$ a suitable quadrangulation
contradicting the minimality of $(S,G,c)$. We mention here that if two pairs of
neighboring edges of $F$ got removed we would end up with $S'$ empty, but this
comes only from the path $P_3$ as an even quadrangulation of the 2-sphere
so it is is not possible.

Next we consider the case that $d$ is separating, but no edge of $F$ got
removed. Now $S'$ is the union of two surfaces having only $x$ as their common
point. The graph $G'$ quadrangulates both
surfaces. As the total number of edges breaking consistency of orientation 
in $G'$ is odd, the subgraph $G''$ of $G'$ quadrangulating one of these
surfaces $S''$ is an odd quadrangulation. Clearly, $S''$ is a
non-orientable surface that has an odd quadrangulation $G''$ and $G''$ has
fewer faces than $G$ but it also inherits a local $3$-coloring. As $(G,S,c)$
is a minimal quadrangulation, the genus of $S''$ must be strictly larger than
that of $S$. But this is impossible as the sum of the Euler characteristics 
of the two surfaces $S'$ consists of is exactly $2$ more than the Euler
characteristic of $S$, and the other surface in $S'$ can contribute at most
$2$ to the sum.

Our remaining last case is when the diagonal $d$ is a non-separating $2$-sided
cycle on $S$. Here $S'$ is a ``pinched surface'', it can be made into a
surface by replacing $x$ with two points. Let $S''$ be the
surface so obtained and let $G''$ be the quadrangulation of $S''$ obtained
in this process. The Euler characteristic of $S''$ is two more than that of $S$
and its quadrangulation $G''$ is odd and has a local $3$-coloring, thus
$(G,S,c)$ is not minimal in this case either. This completes our proof that
$F$ must have four distinct vertices.

Let the vertices along the facial cycle around $F$ be $x$, $y$, $z$ and
$t$. As $c$ is a proper coloring, only the colors of $x$ and $z$, or the colors
of $y$ and $t$ can coincide. In order to prove that these vertices cannot have
exactly $3$ distinct colors assume for a contradiction that $c(x)=c(z)$ and
$c(y)\ne c(t)$. We do as above: we cut $F$ out from $S$ and close $S$ up by
identifying $x$ with $z$ and also the edge $xy$ with the edge $zy$ and the
edge $xt$ with $zt$. This time we obtain a surface $S'$ homeomorphic to
$S$. We also obtain a quadrangulation $G'$ of $S'$ with one fewer faces and
just as above, it must be an odd quadrangulation. As we
identified vertices with equal color, the graph $G'$ inherits a proper
coloring from $c$. We claim it is a local $3$-coloring. This is because both
$x$ and $z$ had $y$ and $t$ in their neighborhoods, so both neighborhoods must
contain only vertices of color $c(y)$ and $c(t)$. Therefore this is also
true for the vertex obtained by identifying $x$ and $z$. The contradiction
with the minimality of $G$ shows that all faces of $G$ must have two or four
distinct colors at their vertices.

Finally, we have to derive a contradiction from the assumption that the
neighborhood of a vertex $x$ in $G$ is monochromatic. Let $F$ be a face
incident to $x$, let $x$, $y$, $z$ and $t$ be the vertices along its facial
walk. As $y$ and $t$ are neighbors
of $x$ we must have $c(y)=c(t)$. By the earlier part of this lemma this
implies that $F$ must have only two colors, so $c(x)=c(z)$. We apply the same
procedure again: cut $F$ out from $S$, and close $S$ up by identifying 
the edge $xy$ with $zy$ and the edge $xt$ with $zt$. We obtain a quadrangulation 
$G'$ of a surface $S'$ with one fewer faces than in $G$. As before, 
$G'$ is an odd quadrangulation, $S'$ is homeomorphic with $S$, and $c$ gives
rise to a local $3$-coloring of $G'$ as all neighbors of the
common image of $x$ and $z$ have colors that appear in the $G$-neighborhood of
$z$. The contradiction with the minimality of $(G,S,c)$ completes the proof of
the lemma. \qed

\bigskip

We call a face of a minimal quadrangulation {\em bichromatic} or {\em
four-chromatic} depending on the number of distinct colors its vertices
receive.

\begin{lem}\label{4ch}
Two four-chromatic faces of a minimal quadrangulation cannot share an edge.
\end{lem}

\proof Let $(G,S,c)$ be a minimal quadrangulation and let the vertices along
the facial cycles of the faces on the edge $e=xy$ be $x$, $y$, $z$, $t$ and
$x$, $y$, $z'$, $t'$. If both of these faces
are four-chromatic, then $c(z)=c(z')$ as otherwise the vertex $y$ had three
different colors in its neighborhood. Similarly, $c(t)=c(t')$. 
Let us obtain $G'$ from $G$ by removing the edge $e$ and inserting a new
edge connecting $z$ and $t'$. Clearly, $G'$ is also an odd quadrangulation of
$S$ with just as many faces as $G$ and $c$ is a local $3$-coloring of
$G'$. Thus $(G',S,c)$ is a minimal quadrangulation. Both faces incident to
the new $zt'$ edge have three distinct colors. This contradicts 
Lemma~\ref{reduction} and proves the present lemma. \qed

\bigskip

We call a vertex of a minimal quadrangulation $(G,S,c)$ {\em regular} if its
degree is $4$. The rest of the vertices of $G$ are called {\em irregular}. We
define the {\em auxiliary graph} $H$ on the vertex set $V(H)=V(G)$ by connecting
diagonally opposite vertices of the bichromatic faces of $G$. Note that in the
auxiliary graph only vertices of equal color are connected.

\begin{lem}\label{degree}
In a minimal quadrangulation each vertex has degree at least four.
In the auxiliary graph
regular vertices of the quadrangulation have degree $2$, irregular vertices
have degree at least $3$. In particular, vertices of degree $5$ in $G$ have
degree $3$ in $H$.
\end{lem}

\proof Let $x$ be a vertex of a minimal quadrangulation $(G,S,c)$. Its degree
$d_G(x)$ in $G$ is the total number of faces incident to $x$, while its
degree $d_H(x)$ in $H$ is the number of bichromatic faces incident to $x$. By
Lemma~\ref{4ch} at least half of the faces incident to $x$ are
bichromatic. Four-chromatic faces correspond to changes in color as we
consider the neighbors of $x$ in their cyclic order. By Lemma~\ref{reduction}
there must be such a change, and therefore at least two four-chromatic
faces. The statements of this lemma follow. \qed

\begin{lem}\label{h-cycle}
If a component $C$ of the auxiliary graph $H$ of a minimal quadrangulation
$(G,S,c)$ is a cycle, then $G$ has at least two vertices of degree
at least\/ $8$ that are adjacent in $G$ to a vertex of $C$.
\end{lem}

\proof Clearly, all vertices 
of $C$ are identically colored and by Lemma~\ref{degree} they are all regular
vertices.

A color is said to {\em match} a vertex $x$ of $C$ if
it is the color $c(z)$ of a vertex $z$ that is diagonally opposite from $x$ in
a four-chromatic face of $G$. Each vertex of $C$ is regular, so each is
incident to exactly two four-chromatic faces, one on either side of $C$.

We claim that the same colors match every vertex of $C$. To see this it is
enough to prove that the same colors match neighboring vertices along $C$. Let
$x$ and $y$ be neighbors along $C$. As they are connected in $H$,
they appear as diagonally opposite vertices of a bichromatic face
$F_1$ of $G$. Consider a color $c(z)$ that matches $x$ with $z$ diagonally
opposite from $x$ in the four-chromatic face $F_2$. As $x$ is a regular vertex
the bichromatic and four-chromatic faces $F_1$ and $F_2$, both incident to
$x$, must share an edge $xt$. Now $ty$ is an edge of $F_1$ and as $y$ is a
regular vertex, the other face $F_3$ incident to this edge must be
four-chromatic. Let $u$ be the vertex diagonally opposite from $y$ in
$F_3$. Clearly, $c(u)$ matches $y$. We finish the proof of the claim by
observing that the neighborhood of $t$ in $G$ contains $x$, $z$ and $u$, so
these vertices cannot have all distinct colors. Since $F_2$ and $F_3$ are
four-chromatic faces, we have $c(z)\ne c(x)$ and $c(u)\ne c(y)=c(x)$,
therefore we must have $c(z)=c(u)$. 

Let us now fix a color $\alpha$ that matches the vertices of $C$ and consider
the coloring $c'$ of the vertices of $G$ given by $c'(x)=\alpha$ for $x$ in
$C$ and $c'(x)=c(x)$ otherwise.

We claim that $c'$ is a proper coloring of $G$. To see this it is enough to
consider a vertex $x$ of $C$ and prove that no neighbor of $x$ in $G$ has color
$\alpha$. As $\alpha$ matches $x$ we have $\alpha=c(z)$ with $z$ diagonally
opposite from $x$ on a four-chromatic face $F$. The two other vertices on $F$
are neighbors of $x$ and have distinct colors neither of which is $\alpha$. As
$c$ is a local $3$-coloring all other neighbors of $x$ must also have one of
these two colors, so none can have the color $\alpha$. 

It is easy to see that $c'$ is not a local $3$-coloring. If it were,
then $(G,S,c')$ would also be a minimal quadrangulation, but any face $F$ that
has a vertex $x$ in $C$ and a vertex $z$ with $c(z)=\alpha$ as diagonally opposite vertices
would have $3$ distinct colors, contradicting Lemma~\ref{reduction}.

We thus have a vertex $w$ that has at least $3$ different colors under $c'$ in
its neighborhood. This is only possible if $w$ has some vertex $x$ of $C$ in
its neighborhood and also some vertex $x'$ of the same color $c(x)=c(x')$ but
outside $C$. Consider the faces around $w$. We saw that four-chromatic faces
correspond to alternations in the color of the neighbors of $w$, so there must
be an even number of those faces. As vertices of $C$ are regular we see that
there are two four-chromatic faces separated by a bichromatic face with two
neighbors of $w$ on $C$. As we must also have $x'$ among its neighbors this
further means that four-chromatic faces incident to $w$ but not to its just
mentioned two neighbors on $C$ must exist, so we have
at least $4$ four-chromatic faces around $w$. By Lemma~\ref{4ch} these faces
must be separated by bichromatic faces, so the degree of $w$ is at least $8$. 

To finish the proof of this lemma we have to establish that the vertex $w$
found 
above is not the only high degree vertex in the neighborhood of $C$. To see
this, notice that $w$ has no neighbor $s$ with $c(s)=\alpha$ as otherwise it
would have no more $c'$-colors in its neighborhood than $c$-colors. But $w$
has a neighbor $x$ in $C$ and one of the two faces incident to the edge $xw$
is a four-chromatic face $F$. If $z$ is diagonally opposite from $x$ in $F$,
then $\beta=c(z)$ matches $x$ and we can define the coloring $c''$ by
re-coloring the vertices of $C$ to $\beta$. With this we find another high
degree vertex as $w$, this one also in the neighborhood of $C$, but not having
color $\beta$ in its neighborhood. \qed

\bigskip
\noindent{\bf Proof of Theorem~\ref{main}(i):}
We shall prove that for a minimal quadrangulation $(G,S,c)$ the genus
$g$ of $S$ is at least $5$. The Euler characteristic of the non-orientable 
surface of genus $g$ is $2-g$. As $G$ is
a quadrangulation, it has half as many faces than edges, so we have
$2-g=|V(G)|+|E(G)|/2-|E(G)|$. Here $2|E(G)|=\sum d_G(x)$, where $d_G(x)$ is
the degree of the vertex $x$ in $G$. We have
$$\sum_{x\in V(G)}(d_G(x)-4)=4(g-2).$$
We call $d_G(x)-4$ the {\em excess} of the vertex $x$. Clearly, regular
vertices have zero excess, irregular vertices have positive excess. Note that
at this point (or rather at Lemma~\ref{degree}) we have re-proved
Theorem~\ref{pplane}, as for $g=1$ the total excess should be $-4$.

We distribute the excess of irregular vertices to the colors. If a vertex has
excess $1$ (i.e., it is of degree $5$), we give this excess to its
color. If the excess of a vertex is $2$ or $3$ we give $2$ of it to its
color. If the excess of a vertex is at least $4$ we still give $2$ to its
color and we give $1$ to each of the two colors in its neighborhood. (Recall
that any vertex has two colors in its neighborhood by Lemma 4.1.)

We have distributed not more, than the total excess of $4(g-2)$. We claim that
each color that is used by the coloring $c$ receives at least $2$ units of the
overall excess. Indeed, the vertices of any color form one or more components
of $H$. If such a component is a cycle, then by
Lemma~\ref{h-cycle} the corresponding color receives $1$ unit of excess from
at least two  distinct high degree neighboring vertices. Now consider a
non-cycle component. By
Lemma~\ref{degree} each vertex in $H$ has degree at least $2$ and degree $3$
means that $1$ unit of excess is given to the color of this vertex, while
degree $4$ or higher means that $2$ units of excess is given to its color. As
no (finite) component can have a single 
degree 3 vertex with all other vertices being of degree $2$, this proves the
claim. 

By Theorem~\ref{am}, $G$ is not colorable by $3$ colors. This implies that $c$
uses at least $5$ colors as any graph that has a local $3$-coloring with $4$
colors can, in fact, be properly colored with $3$ colors (as we have
observed after Theorem~\ref{main}). So the total excess
of $4(g-2)$ is at least $5\cdot 2=10$, and we have $g\ge5$ as claimed. \qed

\section{Constructions}\label{sec:const}

In this section we construct several odd quadrangulations that are locally
$3$-colorable. In particular, we prove Theorems~\ref{main}(ii) and
\ref{5colors}(ii). Our starting points are the following graphs $U(m,r)$, 
defined in \cite{EFHKRS}, characterizing local $r$-colorability.

Let $m\ge r$ be positive integers and $[m]=\{1,\ldots,m\}$. 
The vertex set $V(U(m,r))$ consists of the
pairs $(i,A)$ with $i\in[m]$, $i\notin A\subset[m]$ and $|A|=r-1$.
The vertices $(i,A)$ and $(j,B)$ are adjacent in
$U(m,r)$ if and only if $i\in B$ and $j\in A$. The {\em natural coloring} of
$U(m,r)$ gives the color $i$ to each vertex $(i,A)$. This is a local
$r$-coloring of $U(m,r)$. 

By an elementary result proved in \cite{EFHKRS}, a graph is locally
$r$-colorable using at most $m$ colors if and only if it has a homomorphism to
$U(m,r)$. 

It will be beneficial to distinguish edges of $U(m,3)$ appearing in
triangles. These {\em triangle edges} connect $(i,\{j,k\})$ with
$(j,\{i,k\})$ for some distinct colors $i,j,k\in[m]$.

We let $G_0$ be the subgraph obtained from $U(5,3)$ by removing all triangle
edges. This is an edge-transitive graph on thirty vertices
with sixty edges.

\begin{figure}[htb]
\epsfxsize9cm \centerline{\epsffile{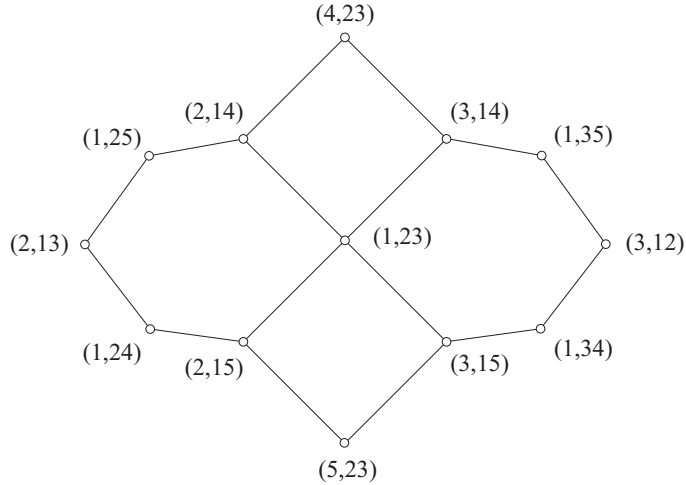}} 
\caption{The faces around the vertex $(1,23)$ in $G_0$}
\label{fig:3}
\end{figure}

We shall describe an embedding of $G_0$ into a surface by listing the faces.
The embedding will have
quadrilateral and hexagonal faces. We take all $4$-cycles in $G_0$ as
quadrilateral faces and we take the $6$-cycles that receive only two distinct
colors at the natural coloring as the hexagonal faces. Notice that each edge
of $G_0$ appears in exactly one quadrilateral and exactly one hexagonal
face. To check that these faces give rise to a surface,
one has to check that the faces form a disk neighborhood around each vertex.
By transitivity of $G_0$, it suffices to verify this for any vertex 
of $G_0$, and we refer to Figure~\ref{fig:3} for details (where we use the 
notation $(i,jk)$ to denote the vertex $(i,\{j,k\})$).

Altogether, we have fifteen quadrilateral faces and ten hexagonal faces.
This makes the Euler characteristic of the resulting surface $S_0$ equal 
to $30-60+25=-5$. Therefore $S_0$ is the non-orientable surface of genus $7$.

\begin{figure}[htb]
\epsfxsize10.2cm \centerline{\epsffile{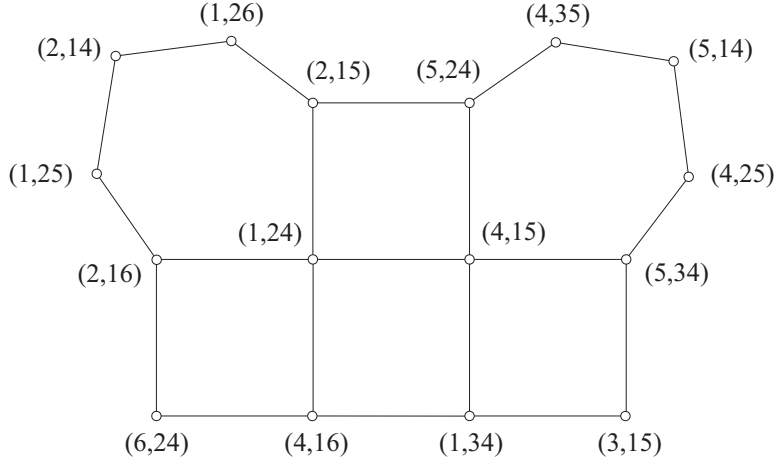}} 
\caption{The faces around vertices $(1,24)$ and $(4,15)$ in $G_1$}
\label{fig:4}
\end{figure}

To obtain a similar example on a surface of smaller genus, we start with the
graph $U(6,3)$. Let $G_1$ be the subgraph of $U(6,3)$ consisting of the
vertices $(i,H)$, with $|H\cap\{1,2,3\}|=1$ (and hence also 
$|H\cap\{4,5,6\}|=1$), and with all non-triangle edges
connecting these vertices. This is a vertex-transitive graph on $36$ vertices
and $72$ edges.

We embed $G_1$ into a surface by listing the resulting faces. As in the case
of $G_0$, we have quadrilateral and hexagonal faces; 
we take all $4$-cycles of $G_1$ as quadrilateral faces and the
$6$-cycles of $G_1$ that receive two colors at the natural coloring as the
hexagonal faces. We have to check again that these faces form a disk
neighborhood around each vertex, and they do (cf.\ Figure~\ref{fig:4}). 
So these faces form a surface $S_1$. We
have $18$ quadrilateral faces that receive four colors at the natural coloring
and nine further quadrilateral faces receiving two colors, and six hexagonal
faces. This makes the Euler characteristic of $S_1$ equal to $36-72+33=-3$.
Therefore $S_1$ is the non-orientable surface of genus $5$.

To obtain quadrangulations of $S_0$ and $S_1$ from the above examples,
we add a main diagonal to every
hexagonal face of $G_0$ and $G_1$. The choice, which of the three main
diagonals to add is arbitrary
for each such face. It is simple to check that the resulting graphs
$G_0'$ and $G_1'$ are odd quadrangulations of $S_0$ and $S_1$, respectively. 
Notice that the
new diagonal edges of $G_0'$ are still edges of the original graph $U(m,3)$
($m=5$ or $6$): they are triangle edges that were removed earlier. This
ensures that the natural coloring locally $3$-colors $G_0'$ with $5$ colors and
locally $3$-colors $G_1'$ with $6$ colors.

We have just given constructions for the first cases of
Theorems~\ref{main}(ii) and \ref{5colors}(ii). To finish the proof of these
results we need to give examples similar to $G_0'$ and
$G_1'$ but quadrangulating higher genus surfaces. For this note that both
quadrangulations $G_0'$ and $G_1'$ have pairs of faces sharing a common edge
and receiving only
two colors. Removing the common edge we get (back) a hexagonal face. We add a
crosscap in the middle of this hexagonal face and quadrangulate the resulting
surface by adding all three diagonals through the crosscap. The surface we
obtain is of genus one higher than our original surface, the resulting
quadrangulation is still odd, and the same coloring still locally $3$-colors
the new graph. Note that the resulting graph will again have neighboring faces
receiving only two colors, so this process can be repeated any number of
times finishing the proof of the existence claims.

We remark that here we increase the genus of the surface with adding new edges
but no new vertices to the quadrangulating graphs. This soon results in
quadrangulating graphs with parallel edges. If one prefers quadrangulating
graphs without parallel edges one can subdivide each edge into three edges and
each face into a grid of nine faces. This has no effect on the parity, but it
gets rid of any parallel edges. Any local $3$-coloring extends to the
subdivided graph using the same set of colors.

We also note that applying the above process exactly once to every hexagonal
face of $G_0$ or $G_1$ we can see that the graph $U(5,3)$ and an induced
subgraph of $U(6,3)$ are themselves odd quadrangulations of some surfaces:
these are the non-orientable surfaces of genus $17$ and $11$,
respectively. \qed

\section{Local $3$-colorings with five colors}\label{sec:5colors}

In this section we turn the construction proving Theorem~\ref{5colors}(ii)
around and use it to prove the impossibility result
Theorem~\ref{5colors}(i). 

Let $G$ be a quadrangulation of a surface and assume each edge of $G$ is
oriented. We call a face of $G$ {\em odd\/} if the edges around it are oriented
with three in one cyclic direction and one in reverse.
The following simple observation shows that odd quadrangulations have an 
odd number of odd faces.

\begin{lem}\label{edgeor}
Let $G$ be a quadrangulation of a surface in which all edges are oriented. 
Then $G$ is an odd quadrangulation if and only if the number of odd faces is odd.
\end{lem}

\proof
Orient the faces of $G$ in an arbitrary manner and consider the facial walks
in this direction. The parity of all face-edge pairs with the
facial walk traversing the edge in reverse direction is the parity of all faces
contributing an odd number -- these are the odd faces -- and also it is the parity
of all edges contributing an odd number -- these are the orientation breaking
edges. \qed
\medskip

\proof(of Theorem~\ref{5colors}(i)).
We consider the graph $G_0$ described in Section~\ref{sec:const} as the graph
embedded in the genus 7 non-orientable surface 
$S_0$. Add all main diagonals in all hexagonal faces to obtain a drawing of
$U(5,3)$ in $S_0$ where each of these newly added  triangular edges cross
two other triangular edges.

Let $G$ be an odd quadrangulation of another surface $S$ and $c$ a local
$3$-coloring of $G$ with the five colors $\{1,2,3,4,5\}$. As we have mentioned
in the previous section (and as proved in \cite{EFHKRS}) $c$ gives rise to a
graph-homomorphism $f:G\to U(5,3)$ such that the natural coloring assigns the
color $c(x)$ to $f(x)$ for every vertex $x$ of $G$. Let $\bar f:S\to S_0$ be a
continuous extension of $f$. First extend $f$ from vertices to the points
along the edges, then realize that the facial walk of any face of $G$ is mapped
trivially to one or two edges of $U(5,3)$ or within a face of $G_0$ and thus
can easily be extended within the same face.

We want to use a simple result that relates the genera of two surfaces and
degree of a mapping between them. As
the degree is usually defined for oriented surfaces we consider the orientable
double cover $\tilde S_0$ of $S_0$ and $\tilde S$ of $S$ with an arbitrary
orientation. Note that the Euler characteristic of the double cover is
twice the Euler characteristic of the base, so the genera of $S$ and $\tilde S$
agree and so do the genera of
$S_0$ and $\tilde S_0$. The map $\bar f:S\to S_0$ 
lifts to $\tilde f:\tilde S\to\tilde S_0$. Let $d$ be the degree of
this map. Consider a quadrilateral face $F$ of $G_0$ and let $n_F$ be the
number of faces of $G$ that $f$ maps to $F$. Notice that the graphs $G$ and
$G_0$ are also lifted to the orientable surfaces and the number of faces
mapped to either one of the faces above $F$ is also $n_F$. To obtain the
degree $d$ of the map $\tilde f$ one has to count these inverse images with
appropriate signs depending on whether $\tilde f$ keeps the orientation of the
face. Thus the
parity of $d$ and $n_F$ must agree. In particular, the parity of $n_F$ is
independent of the choice of the quadrilateral face $F$. (Note that using
hexagonal faces of $G_0$ it would be harder to find the
parity of the degree as some faces of $G$ map to parts of those
hexagonal faces.) We use Kneser's formula \cite{knesersr}, see also
\cite[p.~73]{ZVC}. It states that $g-1\ge|d|(g_0-1)$ for the genera $g>0$ and
$g_0$ of $\tilde S$ and $\tilde S_0$ if a degree $d$ map from $\tilde S$ to
$\tilde S_0$ exists. We use the consequence that $d=0$ whenever $S$
has genus less than $7$.

We claim that $G$ is an odd
quadrangulation if and only if the degree of $\tilde f$ is
odd. As $0$ is even, this claim finishes the proof of Theorem~\ref{5colors}(i).

It is easier to show the claim using Lemma~\ref{edgeor} considered as an
alternative definition of odd quadrangulations.

To obtain an oriented version of $G$ let us orient an edge $xy$ from $x$ to
$y$ if $c(x)<c(y)$. The odd faces will be exactly those with some colors
$a<b<d<e$ around the face in this order. So $G$ is an odd quadrangulation if
and only if $\sum n_F$ is odd, where the summation
is for the quadrilateral faces $F$ of $G_0$ whose verices receive some colors
$a<b<d<e$ in this order in the natural coloring. There are five such faces in
$G_0$. This finishes the proof of the
claim and with it the proof of Theorem~\ref{5colors}(i). \qed

\medskip

Note that a proof similar to the above is not possible for
Theorem~\ref{main}(i) for two reasons. First, $G_1$ is too small a part of
$U(6,3)$, 
there are large parts of $U(6,3)$ that are not represented by $G_1$. Second,
even if $G_1$ was all (or a large part) of $U(6,3)$ this line of thought would
only prove the impossibility of a local 3-coloring of an odd quadrangulation
of small genus surfaces {\em using 6 colors}, whereas Theorem~\ref{main}(i)
states the same for an unlimited number of colors.

\section{Local $4$-colorings of triangulations}
\label{sect:triang}

If $Q$ is a quadrangulation, we define the {\em face subdivision\/} $T(Q)$ 
of $Q$ as the triangulation of the same surface
that is obtained from $Q$ by adding a vertex in each face and joining
it to the four vertices on its boundary. The behavior of the local chromatic
number exhibited in odd quadrangulations also shows up in their face 
subdivisions $T(Q)$.
It has been proved by Hutchinson et al.\ \cite{HRS} that in an odd
quadrangulation 
$Q$ of the projective plane, a vertex coloring with any number of colors gives
rise to  
a four-colored face. This has been generalized to odd quadrangulations of
arbitrary non-orientable surfaces in \cite{AHNNO}; for some further
extensions see also \cite{Mo}. This shows that the face subdivision 
triangulation $T(Q)$ does not admit a local 4-coloring, and its local 
chromatic number is at least 5. Since all added vertices can be colored 
with the same color, we also conclude that $\psi(T(Q))\le \psi(Q)+1$.

\begin{thm}
\label{thm:T(G)}
If\/ $Q$ is an odd quadrangulation of a non-orientable surface and $T(Q)$
is its face subdivision, then the local chromatic number of\/ $T(Q)$ is
at least five.
\end{thm}

There is another family of triangulations of surfaces that exhibits unusual 
chromatic behavior -- a simple local condition forces the chromatic number to 
be at least 5 despite the fact that these graphs seem to be ``almost 3-colorable''. 
Let $T$ be a triangulation of some surface such that all its vertices except two
have even degree, and the two vertices of odd degree are adjacent.
Fisk \cite{Fi} proved that $T$ cannot be 4-colored. 
See \cite[Section 8.4]{MT} for further details.

Here we observe that the above result of Fisk can be extended to local colorings
as well.

\begin{thm}
\label{thm:Fisk}
Let\/ $T$ be a triangulation of some surface such that all its vertices except two
have even degree. If\/ $c$ is a local\/ $4$-coloring of\/ $T$, and\/ $x$ and $y$ 
are the two vertices of odd degree, then $c(x)=c(y)$ and the set of three colors
used on the neighbors of\/ $x$ and of\/ $y$, respectively, is the same.
In particular, if\/ $x$ and $y$ are adjacent, then $T$ has local chromatic 
number at least five.
\end{thm}

\proof
We may assume that $c(V(T))=\{1,2,\dots,m\}$. For every $i,j,k\in\{1,2,\dots,m\}$,
let $T_{ijk}$ be the set of facial triangles whose vertices are colored $i,j,k$,
and let $V_i^{jk}$ be the set of vertices of color $i$ that are incident with a
triangle in $T_{ijk}$.

If $v\in V(T)$, then the set of facial triangles containing $v$ determines
a cycle $C$ (the {\em link} of $v$) passing through all
the neighbors of $v$. Since $c$ is a local 4-coloring, it determines 
a 3-coloring of $C$, which can be viewed as a mapping of $C$ onto the cycle
$C_3$ of length 3. Let $w$ be the winding number of this mapping. Since $C_3$
has an odd number of edges, it follows that 
$$
   w \equiv |C| = \deg(v) \pmod 2.
$$
This simple conclusion implies that the parity of the number of triangles in
$T_{ijk}$ can be expressed as follows:
\begin{equation}
  |T_{ijk}| \equiv \sum_{u\in V_i^{jk}} \deg(u) \equiv 
  \sum_{v\in V_j^{ik}} \deg(v) \equiv 
  \sum_{z\in V_k^{ij}} \deg(z) \pmod 2.
\label{eq:eij}
\end{equation}
There are distinct colors $i,j,k$ such that $x\in V_i^{jk}$ and $y\notin V_j^{ik}$.
The second congruence in (\ref{eq:eij}) can hold only when $y\in V_i^{jk}$
since $x$ and $y$ are the only vertices of odd degree and the sum of degrees 
of vertices in $V_j^{ik}$ is even. This shows that $x$ and $y$ have the same color
and the same set of colors in their neighborhood. 
\qed

\section{Concluding remarks}\label{sec:last}

Quadrangulations can be classified into four types according to the
homology of their odd cycles. Let us describe this refined classification
of even/odd quadrangulations more closely. 
Let $G$ be a graph that is embedded in a surface $S$ such that
all facial walks are of even length. Then it is easy to see that the lengths of 
any two closed walks in $G$ that are homologous (with respect to the 
$\2$-homology $H_1(S,\2)$) have the same parity. The parities of walks in
different homology classes thus determine a homomorphism
$\phi: H_1(S,\2)\to \2$, which is called the {\em cycle parity map} of $G$.
If $S$ is a non-orientable surface of genus $g$, then it is homeomorphic to 
the connected sum of $g$ projective planes $Q_1,\dots,Q_g$, and its homology
group $H_1(S,\2)\cong \2^g$ is generated by 1-sided cycles $\alpha_i$ in $Q_i$,
$i=1,\dots,g$. In particular, $\phi$ can be {\em represented\/} by the $g$-tuple
$(\phi_1,\dots,\phi_g)$, where $\phi_i=1$ if $Q_i$ contains a closed walk in $G$
of odd length, and $\phi_i=0$ otherwise. If $\phi_1=\phi_2=\phi_3=1$ and
$\phi_4=0$, then we can replace the chosen basis $\alpha_1,\dots,\alpha_4$
of the homology group by $\alpha_1' = \alpha_1+\alpha_2+\alpha_3$, 
$\alpha_2' = \alpha_1+\alpha_2+\alpha_4$, $\alpha_3' = \alpha_1+\alpha_3+\alpha_4$, and 
$\alpha_4' = \alpha_2+\alpha_3+\alpha_4$, respectively. 
It is easy to see that $\alpha_1',\dots,\alpha_4',\alpha_5,\dots,\alpha_g$
can be represented by disjoint 1-sided simple closed curves in $S$ that
generate $H_1(S,\2)$. Moreover, since $\phi$ is a homomorphism, it follows that
the parity map representation changes from $(1,1,1,0,\dots)$ to 
$(1,0,0,0,\dots)$ under the new generating set. This shows that there is
a representation of $\phi$ in one (and precisely one) of the following four forms:
\begin{eqnarray*}
 \Phi^0 = (0,0,0,\dots,0), && \Phi^1 = (1,0,0,\dots,0), \\
 \Phi^2 = (1,1,0,\dots,0), && \Phi^3 = (1,1,1,\dots,1).
\end{eqnarray*}
We say that $G$ is of type $\Phi^i$ ($i\in\{0,1,2,3\}$) if its cycle parity map
can be represented by $\Phi^i$.
We refer to \cite{NeNa} or \cite{NaNeOt} for a similar treatment with more
details.

It is clear from the definition that $Q$ is of type $\Phi^3$ if and only if
every 1-sided closed walk has odd length. (It is a corollary of this that every
2-sided closed walk has even length.) Let us observe that the quadrangulation
$G_1'$ of the genus 5 non-orientable surface constructed in Section
\ref{sec:const} is of type $\Phi^3$. 
This can be proved as follows. First we observe that $G_1'$ can be represented
by means of local rotations and the signature (cf.~\cite{MT}). The edges of
negative signature are the following ones:
$$
\begin{tabular}{c c c c}
  (2,35)(3,26)&(2,36)(3,25)&(5,26)(6,35)&(5,36)(6,25) \\
  (3,14)(4,35)&(3,14)(4,36)&(1,26)(6,14)&(1,36)(6,14) \\
  (1,25)(5,14)&(1,35)(5,14)&(2,14)(4,25)&(2,14)(4,26)
\end{tabular}
$$
If we remove all these edges, we get a bipartite spanning subgraph $G_1''$
of $G_1'$, and each removed edge joins two vertices that belong to the same
bipartite class in $G_1''$. This shows that all 1-sided closed walks
(i.e. those that traverse an odd number of edges with negative signature)
in $G_1'$ have odd length. Thus, $G_1'$ is of type $\Phi^3$.

We have shown that odd quadrangulations of non-orientable surfaces of genus
at most four have local chromatic number at least four and that for every
surface of genus at least five, there are examples for which this no longer
holds. Let us observe that a quadrangulation is odd if and only if it is either
of type $\Phi^1$ (for arbitrary genus) or it is of type $\Phi^3$ when the genus
is odd. As shown above, our example $G_1'$ of an odd quadrangulation of 
genus 5 that admits a local 3-coloring is of type $\Phi^3$. 
When producing quadrangulations of higher genera, we can switch to type 
$\Phi^1$ (by repeatedly inserting three diagonals into hexagonal faces as 
done in the proof of Theorem \ref{main}(ii)). But we can also stay within 
the type $\Phi^3$ by taking a quadrangulation of type $\Phi^3$ of odd genus
$g$ that admits a local 3-coloring, such that on two adjacent faces only
two colors are used, and then replace those two faces by making a connected
sum with a 3-colorable quadrangulation (minus an edge $e$) of the Klein bottle 
of type $\Phi^2$ in which the two faces sharing the edge $e$ are 2-colored.
This gives a quadrangulation of genus $g+2$ of type $\Phi^3$
that admits a local 3-coloring, and leaves only one unresolved case --
type $\Phi^1$ on the non-orientable surface of genus 5.

\begin{quest}
\label{q:g5}
Is there a quadrangulation of type $\Phi^1$ of the non-orientable surface of 
genus $5$ that admits a local\/ $3$-coloring?
\end{quest}

\bigskip
\bigskip
\par\noindent
{\bf Acknowledgements.} We thank G\'abor Elek for his help in finding the
papers of Baudisch in the literature and for related conversations. Useful
discussions with L\'aszl\'o Feh\'er and G\'abor Moussong are also gratefully
acknowledged.

\end{document}